\documentclass{article}
\usepackage{graphicx} 
\usepackage{graphicx}      
\usepackage{xcolor}
\usepackage{algorithm}
\usepackage{algpseudocode}
\usepackage{graphicx}
\usepackage{amsfonts}
\usepackage{amsmath,mathtools,amssymb,amsfonts}

\newcommand{\norm}[1]{\left\lVert#1\right\rVert}
\newcommand{\xddots}{%
  \raise 4pt \hbox {.}
  \mkern 6mu
  \raise 1pt \hbox {.}
  \mkern 6mu
  \raise -2pt \hbox {.}
}
\DeclareMathOperator*{\SubjectTo}{Subject\phantom{a}to:}
\DeclareMathOperator*{\Minimize}{Minimize:}
\DeclareMathOperator*{\Maximize}{Maximize:}
\DeclareMathOperator*{\argmin}{arg\,min}

\title{Distributed Model-Predictive Energy Management Strategy for Shipboard Power Systems Considering Battery Degradation}
\author{Satish Vedula$^*$, Seyyed Shaho Alaviani$^{**}$ and \\ Olugbenga Moses Anubi$^*$}

\date{}

\begin{document}

\maketitle
\begin{centering}
$*$ Department of Electrical and Computer Engineering, the Center for Advanced Power Systems, Florida State University 
E-mail: \{svedula\ , oanubi\}@fsu.edu
$**$ Center for Advanced Power Systems (CAPS), Florida State University, and National High Magnetic Field Laboratory (MagLab), Tallahassee, FL 32310, USA (email: salaviani@fsu.edu)
\end{centering}

\section{Abstract}
With the integration of loads such as pulse power loads, a new control challenge is presented in meeting their high ramp rate requirements. Existing onboard generators are ramp rate limited. The inability to meet the load power due to ramp rate limitation may lead to instability. The addition of energy storage elements in addition to the existing generators proves a viable solution in addressing the control challenges presented by high ramp rate loads. A distributed energy management strategy maximizing generator efficiency and minimizing energy storage degradation is developed that facilitates an optimal adaptive power split between generators and energy storage elements. The complex structure of the energy storage degradation model makes it tough for its direct integration into the optimization problem and is not practical for real-time implementation. A degradation heuristic to minimize absolute power extracted from the energy storage elements is proposed as a degradation heuristic measure. The designed strategy is tested through a numerical case study of a consolidated shipboard power system model consisting of a single generator, energy storage element, and load model. The results show the impact of the designed energy management strategy in effectively managing energy storage health.

\textbf{Keywords:} Battery degradation, model predictive control, energy management

\section{Introduction}
Conceptualization of Microgrids (MGs) as an amalgamation of fossil-fuel-driven micro-turbines, renewable energy-based generating systems, and energy storage systems (ESSs) incorporated into a distribution generation environment, driving high demand for the power. MGs are categorized into \emph{grid-connected} and \emph{islanded}, based on the configuration and their mode of operation (\cite{Asanso_2007}). One such example of an islanded operation is Shipboard Power Systems (SPSs). In shipboard power systems, power generation modules (PGMs) and power conversion modules (PCMs) provide the power required by the loads through means of the direct current (DC) distribution system. Modern SPSs are equipped with advanced power loads such as rail guns, electromagnetic radars, and heavy nonlinear and pulsed power loads (PPLs) (\cite{Derry_1}). 

There is a challenge to utilize the existing PGMs, which are ramp-rate limited, to meet the power demanded by the high-ramp-rate loads. Failure to ramp up in appropriate time to meet the load requirements can lead to system imbalances and instability. Furthermore, adding additional generators is not a plausible solution considering the ship hull dimension restrictions. This imposes a need for high-ramp-rate power generation elements to be integrated into SPSs. Energy storage systems such as batteries can support higher ramp rates. Integrating energy storage systems (ESSs) into the existing notional SPS framework addresses the high ramp requirements (\cite{ESRDC_1}). The ESSs can store the energy during normal periods and when needed can provide immediate ramp-up support. Several investigators have addressed the above challenge \emph{without} considering battery degradation (\cite{2017_Vu_2,2021_Vedula,Seenumani}) where they have used dynamic programming (DP) and model predictive control (MPC) to address the aforementioned challenge. MPC is a renowned concept used in control literature whose applications include industry. MPC problem solves an objective function subject to a set of constraints that can be actuator limitations in physical systems or rate constraints (\cite{Wang_Boyd}). The main \emph{advantage} of the MPC over DP is that it enables incorporation of the constraints such as inequality and box constraints capturing the upper and lower limitations and ramp rates, thus accounting for actuator limitations. 

\begin{figure}[t!]
      \centering
      \includegraphics[width=0.70\textwidth]{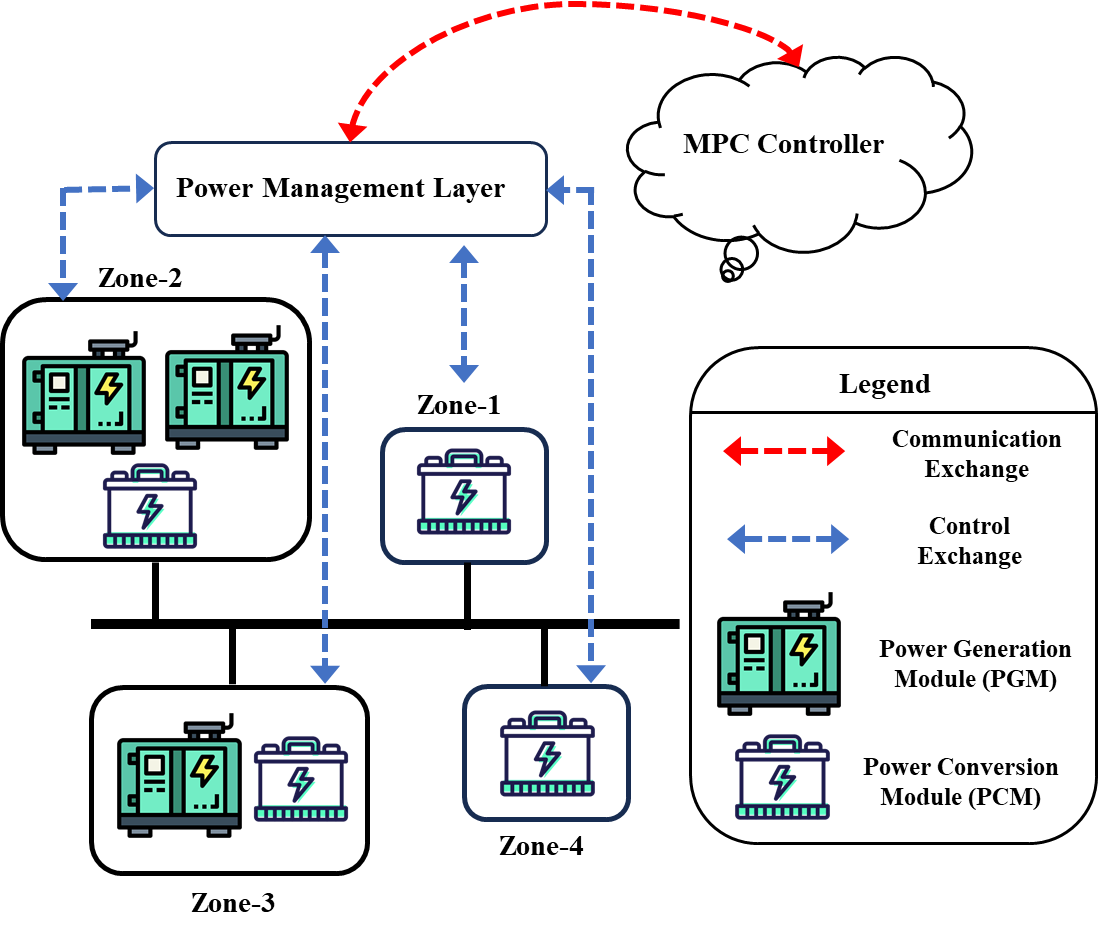}
	 \caption{A Notional Shipboard Power System with the Zonal Architecture and Control Hierarchy.}
     \label{4_Zone_SPS} 
\end{figure}

Nevertheless, batteries degrade faster than generators. This fact imposes another challenge to consider both \emph{ramp rate requirements} and \emph{battery degradation} together. Some investigators have considered this challenge while addressing the EM problem for the MGs (\cite{Hein_2021,Steen_2021,Li_2023,Zhao_2023,Nawaz_2023,Ji_2023,Wang_2022,Satish_2}). In \cite{Hein_2021}, the authors use DP to solve the EM problem while considering the depth of discharge (DoD) as a measure for capturing battery degradation. In \cite{Steen_2021}, an energy market pricing-based EM problem for an MG, considering the battery DoD as a degradation measure, is presented. Market-energy-pricing-based approach for MG with fast Li-ion(standing for Lithium-ion)-based batteries, considering DoD as a measure for capturing the battery degradation, has been presented in \cite{Li_2023}. DoD has been also used as a battery degradation measure to solve the EM problem in \cite{Zhao_2023} where the authors proposed a neural-network-based method to address the EM problem in the MGs. In \cite{Nawaz_2023}, the authors have presented a state-of-charge (SoC)-based EM scheme for MGs; the SoC is regulated around a reference value to minimize battery degradation. Monitoring the battery charge/discharge cycles and using them to determine the battery degradation with an application to the MGs has been investigated in \cite{Ji_2023}. The work of \cite{Hein_2021,Steen_2021,Li_2023,Zhao_2023,Nawaz_2023,Ji_2023} consider battery aging in \emph{different measures} such as DoD, SoC, and charge/discharge based cycle count. In \cite{Wang_2022}, the measures used to minimize the battery usage are \emph{battery power} and \emph{battery SoC}, but the authors in \cite{Wang_2022} did \emph{not} consider the cost associated with the generators and its effect on determining the battery degradation. 

The main contributions of this paper are:
\begin{itemize}
    \item A model predictive energy management strategy is designed to incorporate battery degradation in the form of minimizing its power as a heuristic (measure) to mitigate battery deterioration as shown in Fig.~\ref{4_Zone_SPS}.
    \item We Consider both the maximizing generator efficiency and minimizing battery degradation in the designed model predictive energy management problem along with the component ramp rate limitations to test the performance under the presence of pulsed power loads.
    \item A scalable distributed plug-and-play model predictive energy management strategy is developed and tested on a shipboard power system with the real component ratings provided by the Office of Naval Research (ONR). 
\end{itemize}

The organization of the rest of the paper is as follows: Mathematical notations used throughout the paper are presented in Section \ref{sec:notations}. Section \ref{Sec: Model} presents the physics-based mathematical model of the shipboard power system components and their device-level control. The main results are provided in Section \ref{Sec: Control_Development}. The trade-off between the generator and the energy storage and deductions from the numerical simulation is presented in Section \ref{Sec: Simulation} followed by the concluding remarks in Section \ref{Sec: Conclusion}. 

\section{Notations} \label{sec:notations}
The set of natural numbers, positive real numbers, and real numbers is represented by $\mathbb{N}$, $\mathbb{R}_+$, and $\mathbb{R}$. $\mathbb{L}_{\infty}, \mathbb{L}_2$ denote the signal space of bounded and square-integrable signals. A real matrix with $a$ rows and $b$ columns is denoted as $\mathbb{R}^{a \times b}$. $\mathbb{R}_+$ denotes a positive real number. Scalars are denoted by lowercase alphabets respectively (for example $x \in \mathbb{R}$). The notation $\mathcal{X} \subset \mathbb{R}^n$ represents polytopes of the form $\{\mathbf{x}\in\mathbb{R}^n| A_i\mathbf{x}_i \preceq \mathbf{b}_i, A_e\mathbf{x}_e =\mathbf{b}_e\}$, where $A_i, A_e,\mathbf{b}_i,\mathbf{b}_e$ are the parameters derived from the system specific inclusion constraints. $\underline{\mathbf{0}}$ and $\mathbf{1}$ denote the vector of zeros and ones. For any vector $\mathbf{x} \in \mathbb{R}^n$, $\|\mathbf{x}\|_2 = \sqrt{\mathbf{x}^\top\mathbf{x}}$ and $\|\mathbf{x}\|_1=\sum_{i=1}^{n}|\mathbf{x}_i|$, represents the 2-norm and the 1-norm, respectively (where $|.|$ denotes absolute value). The symbol $\preceq$ denotes the component-wise inequality i.e. $\mathbf{x} \preceq \mathbf{y}$ is equivalent to $\mathbf{x}_i \leq \mathbf{y}_i$ for $i=1,2,\hdots,n$. For a function $f:\mathbb{R}^n \longrightarrow \mathbb{R}^m$, $\nabla f(x)$ denotes the gradient of the function $f$ at $x$. The energy of the signal $x:\mathbb{R}_+ \longrightarrow \mathbb{R}$ is defined as: $$E(x) = \int_0^{\infty}|x(t)|^2dt.$$

\section{Model Development} \label{Sec: Model}
\begin{figure}[t!]
      \centering
      \includegraphics[width=0.65\textwidth]{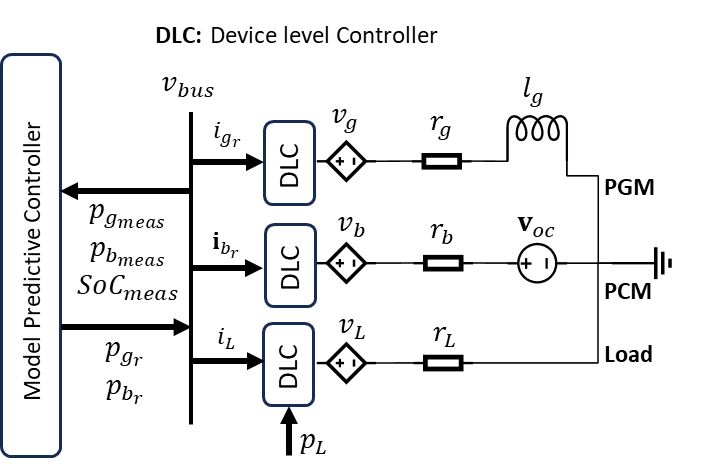}
	 \caption{Mathematical topology of the power exchange in an SPS with a circuit level schematic of PGM, PCM, and Load models.}
     \label{Mathematical_Model} 
\end{figure}

In this section, we present the model development and the device level controller structure for a  Shipboard Power System components \emph{generators/PGMs} and \emph{batteries/PCMs} supplying a common load through a unified DC bus based on Fig.~\ref{Mathematical_Model}.  

The \emph{DC Generator or the PGM} model consists of a current-controlled DC voltage source coupled with a series $RL$ impedance connected to the bus. The assumption is that the bus voltage is already regulated to a set point. Thus, it is considered to be a constant. The dynamical model PGM is given as follows
\begin{equation}\label{gen_dynamics}
 l_{g} \frac{d i_g}{dt} = -r_{g}  i_g(t)+\Delta v(t)
\end{equation}
where $i_g(t) \in \mathbb{R}$ is the PGM current and $\Delta v(t)={v}_{bus}(t)-v_g(t)$, where $v_g(t) \in \mathbb{R}$ is the controllable voltage source and $v_{bus}(t) \in \mathbb{R}$ is the bus voltage to which the PGM is coupled. The generator inductance $l_g \in \mathbb{R}_+$ is in \textsf{Henry}, $r_g \in \mathbb{R}_+$ is the PGM resistance in \textsf{Ohm}. The current injected by the PGM into the bus is dictated by controlling the voltage source.  Given an optimal power profile $p_{g_r} \in \mathbb{R}$ from the energy management control layer, the reference current $i_{g_r} \in \mathbb{R}$ (assumed to be bounded $\in \mathbb{L}_{\infty}$) is generated as $i_{g_r}=p_{g_r}/{v_{bus}}$ and the local controller input $v_g$ is determined via a closed-loop control scheme in which the device level state tracks the reference current i.e. $i_g(t) \to i_{g_r}(t)$. Thus, the DLC is designed to achieve the following objective
$$\int_{0}^{\infty}\bigg(\underbrace{i_g(t)-i_{g_r}(t)}_{\tilde{i}_{g}}\bigg)^2 dt < \infty.$$
Since the design and analysis of such a DLC is not the main objective of this paper, numerous linear (PID) and nonlinear methods can be employed to design it (\cite{Khalil_Book}).

The \emph{Battery Energy Storage System} (BESS) or the \emph{PCM} model consists of multiple ESSs modeled as a voltage-controlled current source. The static model consists of a resistance $r_b \in \mathbb{R}_+$ coupled to the bus $v_{bus} \in \mathbb{R}$ in series with a controllable voltage source $v_b \in \mathbb{R}$ and an open circuit voltage $v_{oc} \in \mathbb{R}$ that make up the battery model. The algebraic equation governing the power exchange in the battery is given as:
\begin{align*}v_{b} &= \frac{v_{bus}^2-p_{b}r_b-v_{bus}v_{oc}}{v_{bus}}, \\
i_b &= \frac{v_{bus}-v_b-v_{oc}}{r_b}.\end{align*}
The state of charge dynamics (SoC) of the battery is given as:
\begin{equation}\label{SoC}
    \frac{d SoC}{dt} = -\frac{1}{Q}i_b(t),
\end{equation}
where $Q \in \mathbb{R}_+$ is the battery capacity in \textsf{ampere-hour}. The discretized form of the SoC dynamics given in (\ref{SoC}) is:
\begin{equation}\label{SoC_Power_Discrete}
SoC_{k+1}=SoC_k-\frac{T_d} {Q}i_{b_k},\end{equation} where $T_d \in \mathbb{R}_+$ is the discretization time-step. Here, $k,k+1$ represents the discrete time intervals.

An Arrhenius equation-based model which uses $Ah$-throughput as $\displaystyle \int_{0}^{t}\left|i_b(\tau)\right|d\tau$ as a metric to evaluate battery state of health (SoH) is used as a PCM degradation metric. $i_b(t) \in \mathbb{R}$ is the current drawn from the individual PCM (positive while discharging and vice-versa). The PCM \textit{capacity loss} formulation is given as follows \cite{SONG2018433}.
\begin{equation}
    Q_L(t) = \zeta_1 e^{\frac{-\zeta_2+TC_{r}}{RT}}\int_{0}^{t}\left|i_b(\tau)\right|d\tau,
\end{equation}
where $T \in \mathbb{R}_+$ is the PCM baseline operating temperature in \textsf{Kelvin}, $C_r$ is the C-rate of the PCM. The capacity loss (\%) is given as follows: $$\Delta Q \% = \frac{Q-Q_{L}(t)}{Q}\times 100,$$
where $Q_{L}$ is the \emph{capacity loss} of the BESS in \textsf{ampere-hour}.

The \textit{Power Load} is modeled as a static resistive load $r_L$ model which consumes the power generated by the PGMs and the PCMs. A variable voltage $v_L \in \mathbb{R}$ enables the current flow in the load based on the active power profile $p_L \in \mathbb{R}$. The assumption is that there is only one load in the model which is consuming all the generated power. Thus, The load model is given as follows:
\begin{align}\label{power_load}
     v_L = \frac{v_{bus}^2-p_Lr_L}{v_{bus}}, \hspace{2mm}
     i_L = \frac{v_{bus}-v_L}{r_L},
\end{align}

Thus, the power flow in the SPS is given by the following algebraic equation:
\begin{equation}\label{PowerFlow}
    \sum_{i=1}^{n_g} \mathbf{p}_{g_i}+\sum_{j=1}^{n_b} \mathbf{p}_{b_j}-p_L = 0,
\end{equation}

where $n_g, n_b \in \mathbb{N}$ are the number of PGMs and PCMs in SPS. Fig.~\ref{Mathematical_Model} shows the hierarchical control structure in the SPS and the power exchange. The power and the SoC measurements are sent to the MPC layer. The MPC computes the optimal values based on the objective and set of constraints and the optimal values act as a reference for DLCs. 

\section{Proposed Energy Management Method}\label{Sec: Control_Development}
This section presents the control development and the derivation of the proposed EM strategy. First, we present the MPC problem formulation in the general form. We make use of the convex analysis and the preliminaries presented in Section \ref{sec:notations} in MPC development. Then, we consider the objective functions of the generator and the battery respectively, and explain the reason behind choosing the particular objectives. Finally, the distributed plug-and-play formulation for the centralized MPC problem is derived. Consider the MPC problem of the form:
\begin{equation}\label{MPC_main}
\begin{aligned}
\Minimize_{\mathbf{p}_{g_i},\mathbf{p}_{b_j}} \quad & \sum_{{i}=1}^{n_g}{C}_{g_{i}}(\mathbf{p}_{{g}_{i}})+\sum_{{j}=1}^{n_b}{C}_{b_{j}}(\mathbf{p}_{{b}_{j}})\\
\SubjectTo \quad & \sum_{{i}=1}^{n_g}\mathbf{p}_{{g}_{i}}+\sum_{{j}=1}^{n_b}\mathbf{p}_{{b}_{j}} - p_f\mathbf{1} = 0, \hspace{1mm} \forall  k,\\
& \mathbf{p}_{{g}_{{i}}} \in \mathcal{X}_g, \hspace{1mm} \forall k, \\& \mathbf{p}_{{b}_{j}}  \in \mathcal{X}_b, \hspace{1mm}\forall k,
\end{aligned}
\end{equation}
where $k \triangleq \left[\begin{array}{cccc} 1&2&\hdots&h\end{array}\right]^{\top} \in \mathbb{R}^h$ is the length of the prediction horizon, $\mathbf{p}_{g_i} \triangleq \left[\begin{array}{cccc}{\mathbf{p}_{g_i}}_k&{\mathbf{p}_{g_i}}_{k+1}&\hdots&{\mathbf{p}_{g_i}}_{k+h-1}\end{array}\right]^{\top}\in\mathbb{R}^h$ is the power profile for generator $i$ over the prediction horizon of length $h$, \newline $\textbf{p}_{b_j} \triangleq \left[\begin{array}{cccc}{\mathbf{p}_{b_j}}_k&{\mathbf{p}_{b_j}}_{k+1}&\hdots&{\mathbf{p}_{b_j}}_{k+h-1}\end{array}\right]^{\top}\in\mathbb{R}^h$ is the power profile for battery $j$ over the prediction horizon of length $h$ and ${p}_f\mathbf{1} \in\mathbb{R}^h$ is the desired total power held fixed over the prediction horizon. The generator cost ${C}_{g_i}:\mathbb{R}^{{h}}\longrightarrow\mathbb{R}_+$ could be associated with an efficiency map or operating at a desired rated power. The battery cost ${C}_{b_j}:\mathbb{R}^{{h}}\longrightarrow\mathbb{R}_+$ could be the cost associated with battery health monitoring and degradation management. ${h} \in \mathbb{N}$ represents the prediction horizon. $\mathcal{X}_g \subset \mathbb{R}^{h}$ and $\mathcal{X}_b \subset \mathbb{R}^h$ represents the set of inclusion constraints each of the PGM and PCM summarizing attributes such as upper and lower power limits, ramp rate limits, and relevant system dynamics. The optimization variables are optimized for the entire length of the horizon and the first value of the sequence is supplied/applied as a power reference (control input) to the low-level (DLC). The PGM and the PCM costs are considered so that the objective is to maintain the PGM around a rated value ($\mathbf{p}_g^r$) which is known and to minimize the PCM power. Based on the optimization problem presented in (\ref{MPC_main}), consider the following MPC problem with the respective cost functions for the PGM and the PCM as follows:
\begin{equation} \label{Full_MPC}
\begin{aligned}
\Minimize_{\mathbf{p}_{g},\mathbf{p}_{b},\mathbf{SoC}_b} \quad & \sum_{i=1}^{n_g}\frac{\beta_i}{2}\norm{\mathbf{p}_{g_{i}}-\mathbf{p}^{r}_{{g_i}}}_2^2+\sum_{j=1}^{n_b}\frac{\gamma_j}{2}\norm{\mathbf{p}_{b_{j}}}_2^2\\
\SubjectTo \quad & \sum_{i=1}^{n_g}\mathbf{p}_{g_{i}}+\sum_{j=1}^{n_b}\mathbf{p}_{b_{j}} = p_f\mathbf{1}, \forall k, \\
\quad & \mathbf{SoC}_{b_{jk+1}} = \mathbf{SoC}_{b_{jk}}-\frac{T_d}{Qv_{bus}} \mathbf{p}_{b_{jk}}, \forall k, \\ \quad & \underline{\mathbf{p}}_{g} \preceq \mathbf{p}_{g_{i}} \preceq \overline{\mathbf{p}}_{g}, \forall k,\\
& \left|\mathbf{p}_{g_{ik}}-\mathbf{p}_{g_{ik-1}} \right| \preceq {r}_{g}\mathbf{1}, \forall k \\ \quad & \underline{\mathbf{p}}_{b} \preceq \mathbf{p}_{b_{j}} \preceq \overline{\mathbf{p}}_{b}, \forall k,\\
& \left|\mathbf{p}_{b_{jk}}-\mathbf{p}_{b_{jk-1}}\right| \preceq {r}_{b} \mathbf{1}, \forall k, \\
& \underline{\mathbf{SoC}}_b \preceq \mathbf{SoC}_{b_{jk}} \preceq \overline{\mathbf{SoC}}_b, \forall k,
\end{aligned} 
\end{equation}
where $k \triangleq \left[\begin{array}{cccc} 1&2&\hdots&h\end{array}\right]^{\top} \in \mathbb{R}^h$, $r_g \in \mathbb{R}_+$ and $r_b \in \mathbb{R}_+$ are the ramp rates of the generator and the battery respectively. $\underline{\mathbf{p}}_{g}$ and $\overline{\mathbf{p}}_{g}$ are the lower and the upper power limitations on the generator, respectively. $\underline{\mathbf{p}}_{b}$ and $\overline{\mathbf{p}}_{b}$ are the lower and the upper power limitations on the battery. $\underline{\mathbf{SoC}}_b$ and $\overline{\mathbf{SoC}}_b$ represent the lower and the upper limitations on the state of charge of the battery. $\beta_i \in \mathbb{R}_+^{n_g}$ are the tunable parameters penalizing the PGMs deviation from a known efficient operation point. $\gamma_j \in \mathbb{R}_+^{n_b}$ are the tunable parameters penalizing the battery power. 

\subsection{Distributed Formulation}
The optimization problem in  (\ref{MPC_main}) is broken down into the individual optimization problems solved at the respective PGMs and the PCMs coupled via an aggregator.
\begin{align}\nonumber
\mathcal{L}(\textbf{p}_{{g}_{i}},\textbf{p}_{{b}_{j}},\boldsymbol{\lambda}) &= \sum_{{i}=1}^{n_g}{C}_{{g_i}}(\mathbf{p}_{{g}_{i}})+\sum_{{j}=1}^{n_b}{C}_{{b_j}}(\mathbf{p}_{{b}_{j}})+\mathcal{I_X}_g+\mathcal{I_X}_b\\\nonumber
&\hspace{10mm}+\boldsymbol{\lambda}^\top\bigg(\sum_{{i}=1}^{n_g}\mathbf{p}_{{g}_{i}}+\sum_{{j}=1}^{n_b}\mathbf{p}_{{b}_{j}} - {p}_{f}\mathbf{1}\bigg)\\\nonumber
&= \mathcal{I_X}_g+\mathcal{I_X}_b+ \sum_{{i}=1}^{n_g}\bigg({C}_{{g_i}}(\mathbf{p}_{{g}_{i}})+\boldsymbol{\lambda}^\top\mathbf{p}_{{g}_{i}}\bigg) \\\nonumber
&\hspace{5mm}+\sum_{{j}=1}^{n_b}\bigg({C}_{{b_j}}(\mathbf{p}_{{b}_{j}})+ \boldsymbol{\lambda}^\top\mathbf{p}_{{b}_{j}}\bigg)\label{Lagrangian_MPC}
-p_f\boldsymbol{\lambda}^\top\textbf{1},
\end{align}
where $\boldsymbol{\lambda} \in \mathbb{R}^{h}$ is the dual variable associated with the main problem in (\ref{MPC_main}) and $\mathcal{I_X}_g$ and $\mathcal{I_X}_b$ are the indicator functions for the inclusion constraints. Given $\boldsymbol{\lambda}$, let $\mathbf{p}_{g_i}^*$ and $\mathbf{p}_{b_j}^*$ to be the solutions to the optimization problem in (\ref{MPC_main}) (\cite{boyd_vandenberghe_2004}) 
\begin{subequations}\label{Iteration_Updates_1}
\begin{align}
    \textbf{p}^*_{g_i} (\boldsymbol{\lambda}) = \argmin_{{\mathbf{p}_{{g}_{i}}} \in \mathcal{X}_g} \bigg\{\text{C}_{{g_i}}(\mathbf{p}_{{g}_{i}})+\boldsymbol{\lambda}^{\top}\textbf{p}_{{g}_{i}}\bigg\}, \\
    \textbf{p}^*_{b_j}(\boldsymbol{\lambda}) = \argmin_{{\mathbf{p}_{{b}_{j}}} \in \mathcal{X}_b} \bigg\{\text{C}_{{b_j}}(\mathbf{p}_{{b}_{j}})+\boldsymbol{\lambda}^{\top}\textbf{p}_{{b}_{i}}\bigg\}.
\end{align}
\end{subequations}
Thus, the dual problem is given as:
$$\Maximize \quad \mathcal{L}(\textbf{p}^*_{g_i},\textbf{p}^*_{b_j},\boldsymbol{\lambda}) \triangleq d(\boldsymbol{\lambda}), $$
which is solved using the gradient ascent algorithm
\begin{align}\nonumber
    \boldsymbol{\lambda}_{t+1} &= \boldsymbol{\lambda}_t +  \alpha \nabla d(\boldsymbol{\lambda}_t)\\\label{eqn:dual_ascent1}
                               &= \boldsymbol{\lambda}_t + \alpha \left(\sum\limits_{i=1}^{n_g}\mathbf{p}_{g_i}^*+\sum\limits_{j=1}^{n_b}\mathbf{p}_{b_j}^* - p_f\mathbf{1}\right).
\end{align}
for a step size $\alpha > 0$ and $t$ is the \emph{iteration counter}. Since the minimizing power profiles $\textbf{p}^*_{g_i}$ and $\textbf{p}^*_{b_j}$ depends on $\boldsymbol{\lambda}$, the dual ascent step in \eqref{eqn:dual_ascent1} becomes a recursion and challenging to evaluate. Consequently, we incorporate the primal problems in \eqref{Iteration_Updates} within a sequential minimization-maximization scheme \cite{Boyd_ADMM} as follows:
\begin{subequations}\label{Iteration_Updates}
\begin{align}
    \mathbf{p}_{g_i}^{t+1} &= \argmin_{{\mathbf{p}_{{g}_{i}}} \in \mathcal{X}_g} \bigg\{{C}_{{i}}(\mathbf{p}_{{g}_{i}})+\boldsymbol{\lambda}^{\top}_t\mathbf{p}_{{g}_{i}}\bigg\}, \\
    \mathbf{p}_{b_j}^{t+1} &= \argmin_{{\mathbf{p}_{{b}_{j}}} \in \mathcal{X}_b} \bigg\{{C}_{{j}}(\mathbf{p}_{{b}_{j}})+\boldsymbol{\lambda}^{\top}_t\mathbf{p}_{{b}_{j}}\bigg\},\\
    \boldsymbol{\lambda}_{t+1} &= \boldsymbol{\lambda}_{t}+\alpha \bigg(\sum_{{i}=1}^{n_g}\mathbf{p}^{t+1}_{{g}_{i}}+\sum_{{j}=1}^{n_b}\mathbf{p}^{t+1}_{{b}_{j}} - p_f\mathbf{1}\bigg). 
\end{align}
\end{subequations}

Next, more details on the resulting nodal problems in the above scheme with the choice of the cost functions and the inclusion constraints are given.
\subsection{PGM Node Optimization Problem}
From (\ref{Iteration_Updates}a) the optimization problem at the $i^{th}$ PGM node assuming the cost function $C_i(\mathbf{p}_{g_i})$ to be a quadratic function is given in (\ref{PGM_Node}). The objective of the considered cost is to maintain the PGM power at around a given rated value or a set point value at every instant of the horizon given by $\mathbf{p}^r_{g_i}$. This rated value is the operating point for that generator as prescribed by the manufacturer. The optimization problem is given as follows:
\begin{equation} \label{PGM_Node}
\begin{aligned}
\Minimize_{\mathbf{p}_{{g}_{i}}} \quad & \frac{\beta_i}{2}\norm{\mathbf{p}_{{g}_{i}}-\mathbf{p}^r_{{g_i}}}_2^2+\boldsymbol{\lambda}^{\top}_t  \mathbf{p}_{{g}_{{i}}}\\
\SubjectTo \quad & \underline{\mathbf{p}}_{g_i} \preceq \mathbf{p}_{{g_i}_k} \preceq \overline{\mathbf{p}}_{g_i}, \hspace{1mm} \forall k\\
&\left|\mathbf{p}_{{g_i}_k}-\mathbf{p}_{{{g}_{{i}}}_{k-1}}\right| \preceq {r}_{g}\mathbf{1}, \hspace{1mm} \forall k,
\end{aligned} 
\end{equation}
where ${r}_{g} \in \mathbb{R}_+$ is the ramp-rate limitation associated with generator. $\overline{\mathbf{p}}_{{{{g}_{i}}}} \in \mathbb{R}_+$ and $\underline{\mathbf{p}}_{{{{g}_{i}}}} \in  \mathbb{R}_+$ are the upper and lower limits on generator power respectively.

\begin{figure}[t!] 
	\centering
	\includegraphics[width = 0.75\textwidth]{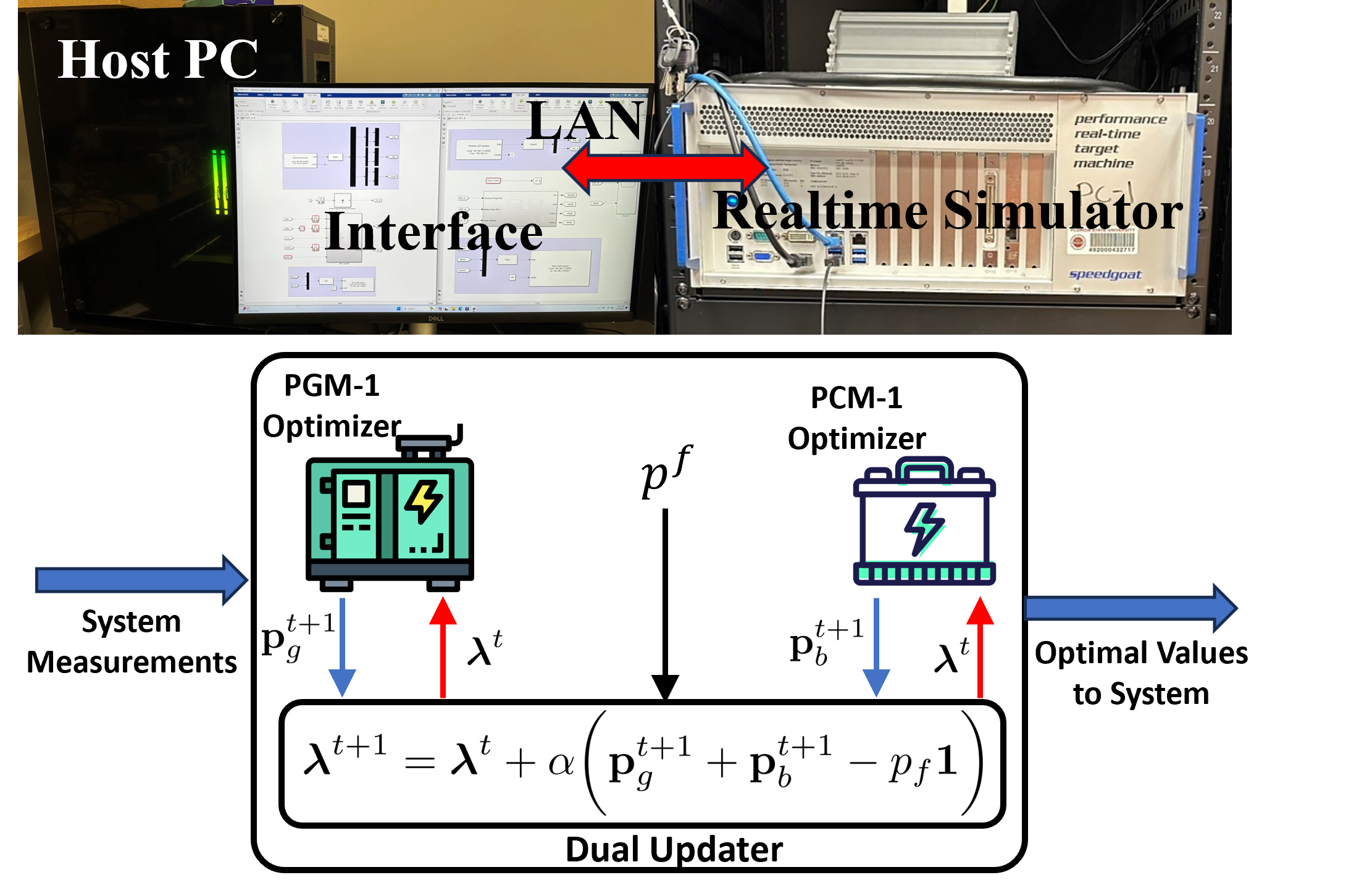} 
	\caption{Distributed Dual Ascent Algorithm Implementation Schematic for Numerical Simulation.
    }
	\label{DA_Algorithm}
\end{figure}

\subsection{PCM Node Optimization problem}
From (\ref{Iteration_Updates}b) the optimization problem at each battery/PCM node is given in (\ref{PCM_Node}). The discretized SoC dynamics in (\ref{SoC_Power_Discrete}) is introduced as an equality constraint. Since capturing the PCM degradation is a long process, the heuristic considered in this work is to minimize the PCM power, which translates into minimizing the PCM degradation. Thus, the cost is introduced as the weighted 2-norm of the PCM power minimization. The goal is to minimize the usage of batteries as much as possible. Thus, the optimization problem is as follows:
\begin{equation} \label{PCM_Node}
\begin{aligned}
\Minimize_{\mathbf{p}_{{b_j}},\mathbf{SoC}_{b_j}} \quad & \frac{\gamma_j}{2}\norm{\mathbf{p}_{b_j}}_2^2+\boldsymbol{\lambda}^{\top}_t  \mathbf{p}_{b_j}\\
\SubjectTo \quad &  \mathbf{SoC}_{b_{jk+1}} = \mathbf{SoC}_{b_{jk}}-\kappa \mathbf{p}_{b_{jk}}, \hspace{1mm} \forall k, \\ & \underline{\mathbf{p}}_{b_j} \preceq \mathbf{p}_{{b_j}_k} \preceq \overline{\mathbf{p}}_{b_j}, \hspace{1mm} \forall k,\\
& \underline{\mathbf{0}} \preceq \mathbf{SoC}_{{b_j}_k} \preceq \mathbf{1}, \hspace{1mm} \forall k, \\
& \left|\mathbf{p}_{{b_j}_k}-\mathbf{p}_{{b_{j}}_{k-1}}\right| \preceq {r}_{b}\mathbf{1}, \hspace{1mm}, \forall k,
\end{aligned} 
\end{equation}
where, $\displaystyle \kappa = {T_d}\slash{Q v_{bus}}$, $r_b \in \mathbb{R}_+ $  is the ramp-rate limitation associated with battery, $\overline{\mathbf{p}}_{b_j} \in \mathbb{R}_+$ and $\underline{\mathbf{p}}_{b_j} \in \mathbb{R}_-$ are the upper and lower limits on battery power. This indicates the discharging and the charging mode of the PCM. The solution to the optimization problems in (\ref{PGM_Node}) and (\ref{PCM_Node}) are used in calculating the dual variable update. The schematic for this aggregation and distribution setup is provided in Fig.~\ref{DA_Algorithm}. All the nodes present in the system send their current iterate value to the dual updater node and the information of the previous iterate of the dual variable $\boldsymbol{\lambda}_t$ is sent to all the nodes.  Thus, the dual updater, which acts as a common scheduler to all nodes implements the update law given in (\ref{Iteration_Updates}c). The convergence of the proposed algorithm under the assumed objectives and the constraint sets which are polytopes is well established in the literature (\cite{boyd_vandenberghe_2004}; \cite{Boyd_ADMM}).

\section{Numerical Simulation}\label{Sec: Simulation}
The designed model predictive energy management strategy is tested on a single PGM, PCM, and Load scenario to fully study and analyze the effect of the weights $\beta$ and $\gamma$ on the total energy delivered by the PCM, PGM, and its effect on the PCM capacity loss \%. The simulation environment used in this setup is based on MATLAB-Simulink on a real-time simulator. The simulation is run on a desktop with the following configurations Digital Storm Intel Core i9-13900K (5.7GHz Turbo) with 64GB RAM acting as a Host PC and a real-time target with Intel Core 3.6GHz 8-core and 32GB RAM (shown in Fig.~\ref{DA_Algorithm}). The optimization problem presented in (\ref{Full_MPC}) with $n_g = 1$ and $n_b=1$ is implemented to study the trade-off characteristics. The parameters used are based on the documentation and the notional SPS component sizing provided by the Office of Naval Research (ONRs) ESRDC (\cite{ESRDC_1}). For all the results presented, the prediction horizon is fixed at ($h=5$\textsf{seconds}). The simulation was run at the fixed step of $10^{-3}$\textsf{seconds}. The rate transition or the communication (measurement exchange) delay between the model and the MPC EM layer was set at $1$\textsf{second} to facilitate optimization convergence. 

\begin{figure}[b!]
      \centering
      \includegraphics[width=0.75\textwidth]{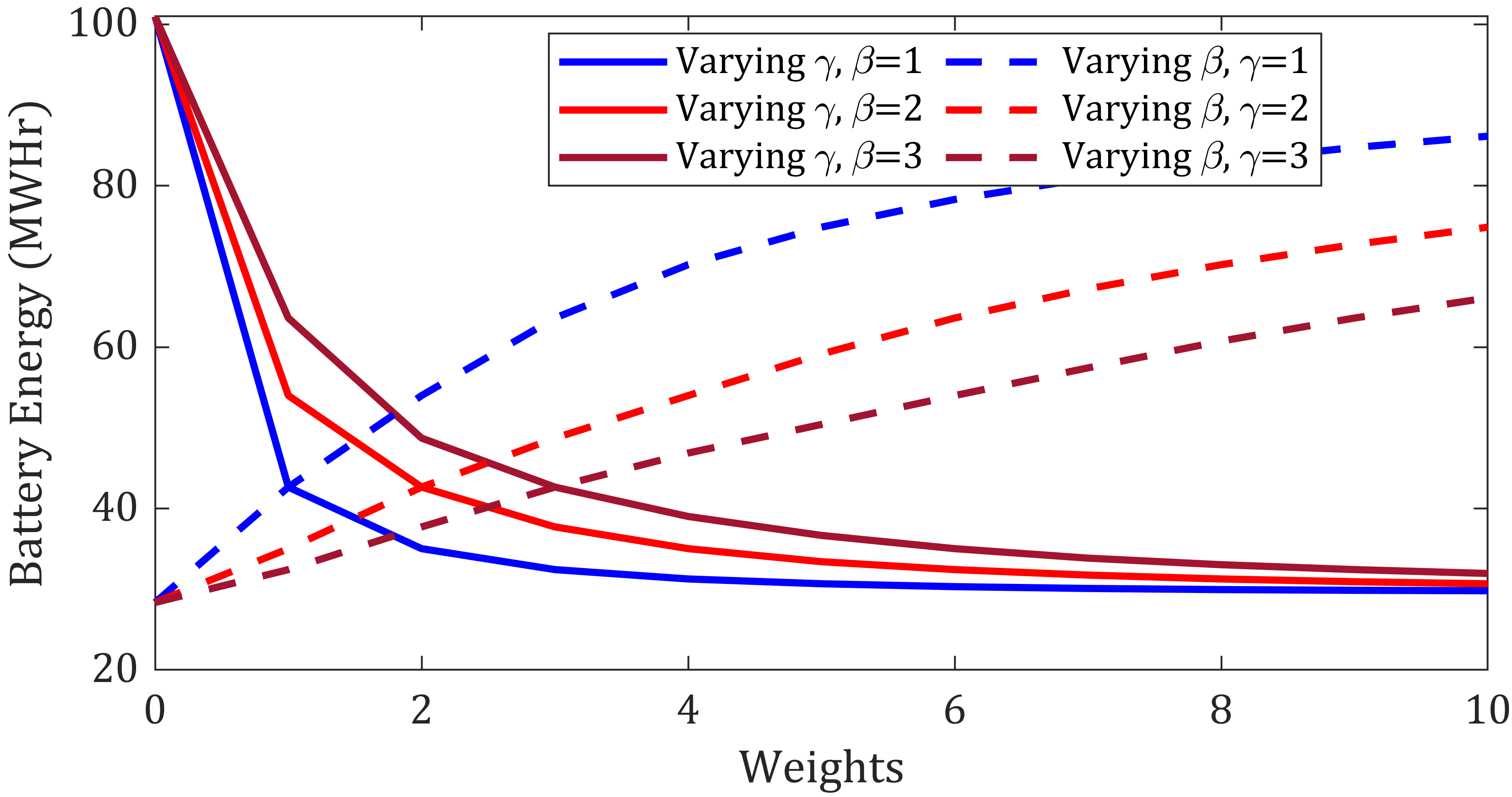}
	 \caption{The effect of different weighting of $\gamma$ and $\beta$ versus the battery energy.}
     \label{P_batt vs Weights} 
\end{figure}

Since the core focus of the work is on determining the effect of the weights $\beta$ and $\gamma$ on the generator and the battery usage and finding a trade-off between them, the results focus on the trade-off part. First, the developed EM algorithm is tested on a \emph{single PGM, PCM, and Load} model. Fig.~\ref{P_batt vs Weights} shows the total energy the battery utilized in \textsf{MWHr} on the y-axis versus the weight on the x-axis. The \emph{solid lines} in the figure present the effect of the battery penalty parameter $\gamma$ tuned from 0-10 against the generator penalty parameter $\beta$. As the $\gamma$ increases the total energy utilized by the battery decreases and can be seen around 30\textsf{MWHr}. This trend shifts upward for the penalty tuning of $\gamma$ from 0-10 but for the increased value of $\beta$. This shift is attributed to the power balance constraint. Similarly \emph{dotted lines} depict the effect of the generator penalty parameter $\beta$ tuned from 0-10 against the battery penalty $\gamma$. It can be seen that as the value of the $\beta$ increases, the total energy utilized by the battery increases since the emphasis is on maintaining the generator around the rated operating point. 

\begin{figure}[h!]
      \centering
      \includegraphics[width=0.75\textwidth]{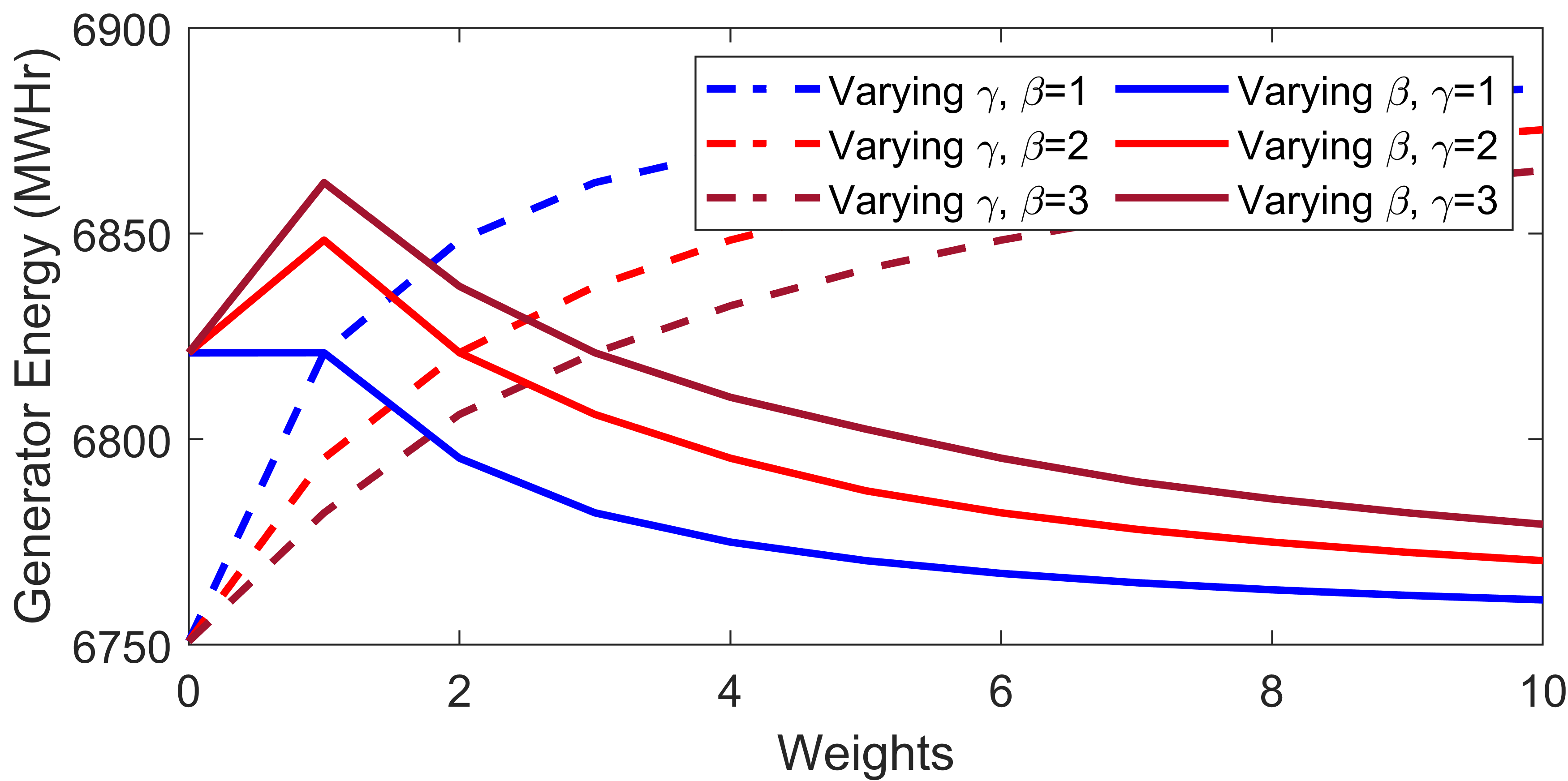}
	 \caption{The effect of different weighting $\gamma$ and $\beta$ versus the generator energy.}
  \label{Gen_Power vs Weights} 
\end{figure}

\begin{figure}[h!]
      \centering
      \includegraphics[width=0.75\textwidth]{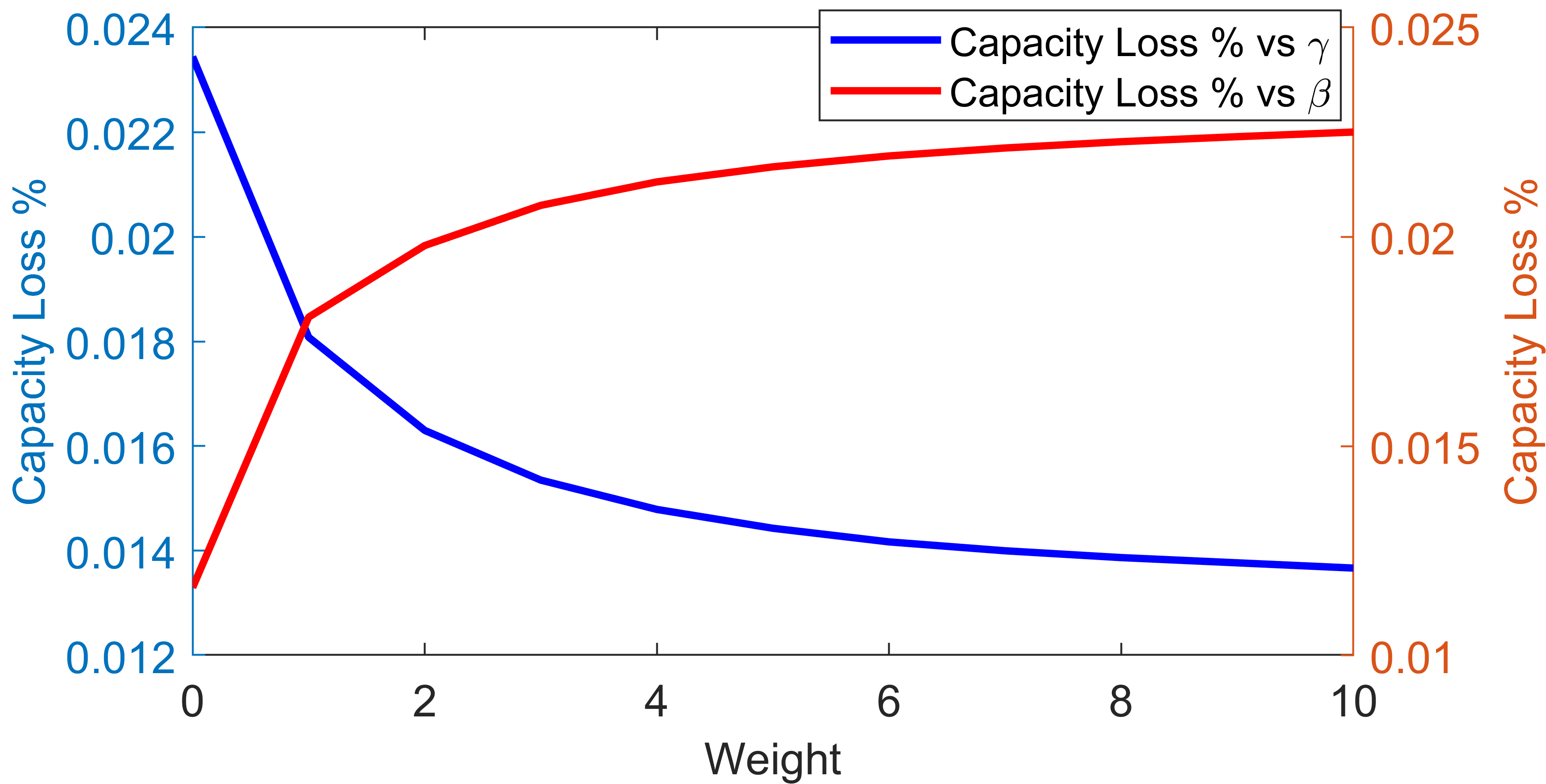}
	 \caption{Battery Capacity Loss $\%$ vs the Weight $\gamma$.}
     \label{Cap_loss vs Beta} 
\end{figure}

In Fig.~\ref{Gen_Power vs Weights} an experiment similar to the previous one is performed, and the energy of the generator is captured versus the weights $\gamma$ and $\beta$. The trends observed for the generator are opposite to the ones observed for the battery. Mathematically this is what is expected to maintain the power supply and demand equality constraint. The \emph{solid lines} depicts the effect of the weight $\gamma$ tuned from 0-10 while keeping $\beta=1$ and then at $\beta=2$ and $3$. When $\gamma=0$, the objective functions goal is minimizing $\frac{\beta}{2}\norm{\mathbf{p}_g-\mathbf{p}_g^r}_2^2$. 
It can be seen that in this case, the energy of the generator is around 6830\textsf{MWHr}. But as the value of the weight $\gamma$ goes to 10, the objective shifts to minimizing the battery usage, and thus, the generator energy goes up to maintain the power balance. Similarly, the \emph{dotted lines} show the effect of weight $\beta$ swept from 0-10 on the generator energy when the weight $\gamma=1$, $\gamma=2$, and $\gamma=3$. Fig.~\ref{Cap_loss vs Beta} shows the effect of the weight $\gamma$ incremented from 0-10 while keeping the weight $\beta=1$ on the battery capacity loss \%. It can be seen that, as the emphasis on the weight $\gamma$ increases the capacity loss percent decreases. In summary, the total capacity loss decreases with $\gamma$ but increases with $\beta$. From the simulated results, we can conclude that a designer can tune the weights $\beta$ and $\gamma$ to obtain the optimal performance of the generator and the battery.

\section{Conclusion} \label{Sec: Conclusion}

A distributed MPC-based energy management strategy considering energy storage degradation is presented for a shipboard power system. The generator, battery ramp rate limitations, generator rated conditions, the pulsed power load conditions are considered. Absolute power extracted from, the battery as a model-based battery degradation heuristic capturing the battery usage is proposed and used in the optimization. Finally, a numerical case study shows the trade-off in the power-sharing and the energy exchange for different optimization weights is presented.

\bibliographystyle{IEEEtran}
\bibliography{myreferences}

\begin{thebibliography}{10}
\providecommand{\url}[1]{#1}
\csname url@samestyle\endcsname
\providecommand{\newblock}{\relax}
\providecommand{\bibinfo}[2]{#2}
\providecommand{\BIBentrySTDinterwordspacing}{\spaceskip=0pt\relax}
\providecommand{\BIBentryALTinterwordstretchfactor}{4}
\providecommand{\BIBentryALTinterwordspacing}{\spaceskip=\fontdimen2\font plus
\BIBentryALTinterwordstretchfactor\fontdimen3\font minus \fontdimen4\font\relax}
\providecommand{\BIBforeignlanguage}[2]{{%
\expandafter\ifx\csname l@#1\endcsname\relax
\typeout{** WARNING: IEEEtran.bst: No hyphenation pattern has been}%
\typeout{** loaded for the language `#1'. Using the pattern for}%
\typeout{** the default language instead.}%
\else
\language=\csname l@#1\endcsname
\fi
#2}}
\providecommand{\BIBdecl}{\relax}
\BIBdecl

\bibitem{Asanso_2007}
N.~Hatziargyriou, H.~Asano, R.~Iravani, and C.~Marnay, ``Microgrids,'' \emph{Power and Energy Magazine, IEEE}, vol.~5, pp. 78 -- 94, 08 2007.

\bibitem{Derry_1}
N.~Doerry and J.~Amy, ``Mvdc shipboard power system considerations for electromagnetic railguns,'' 09 2015.

\bibitem{ESRDC_1}
E.~Team, ``Model description document: Notional four zone mvdc shipboard power system model,'' \emph{ESRDC Website, www.esrdc.com}, 2020.

\bibitem{2017_Vu_2}
T.~V. Vu, D.~Gonsoulin, D.~Perkins, B.~Papari, H.~Vahedi, and C.~S. Edrington, ``Distributed control implementation for zonal mvdc ship power systems,'' in \emph{2017 IEEE Electric Ship Technologies Symposium (ESTS)}, 2017, pp. 539--543.

\bibitem{2021_Vedula}
M.~M. Bijaieh, S.~Vedula, and O.~M. Anubi, ``Model and load predictive control for design and energy management of shipboard power systems,'' in \emph{2021 IEEE Conference on Control Technology and Applications (CCTA)}, 2021, pp. 607--612.

\bibitem{Seenumani}
G.~Seenumani, J.~Sun, and H.~Peng, ``{Exploiting Time Scale Separation for Efficient Real-Time Optimization of Integrated Shipboard Power Systems},'' ser. Dynamic Systems and Control Conference, vol. ASME 2008 Dynamic Systems and Control Conference, Parts A and B, 10 2008, pp. 545--552.

\bibitem{Wang_Boyd}
Y.~Wang and S.~Boyd, ``Fast model predictive control using online optimization,'' \emph{IEEE Transactions on Control Systems Technology}, vol.~18, no.~2, pp. 267--278, 2010.

\bibitem{Hein_2021}
K.~Hein, Y.~Xu, G.~Wilson, and A.~K. Gupta, ``Coordinated optimal voyage planning and energy management of all-electric ship with hybrid energy storage system,'' \emph{IEEE Transactions on Power Systems}, vol.~36, no.~3, pp. 2355--2365, 2021.

\bibitem{Steen_2021}
K.~Antoniadou-Plytaria, D.~Steen, L.~A. Tuan, O.~Carlson, and M.~A. Fotouhi~Ghazvini, ``Market-based energy management model of a building microgrid considering battery degradation,'' \emph{IEEE Transactions on Smart Grid}, vol.~12, no.~2, pp. 1794--1804, 2021.

\bibitem{Li_2023}
J.~Li, X.~Xu, Y.~Wang, R.~Chen, and C.~Liu, ``Bi-level optimizing model for microgrids with fast lithium battery energy storage considering degradation effect,'' \emph{IEEE Access}, vol.~11, pp. 34\,643--34\,658, 2023.

\bibitem{Zhao_2023}
C.~Zhao and X.~Li, ``Microgrid optimal energy scheduling considering neural network based battery degradation,'' \emph{IEEE Transactions on Power Systems}, vol.~39, no.~1, pp. 1594--1606, 2024.

\bibitem{Nawaz_2023}
A.~Nawaz, J.~Wu, J.~Ye, Y.~Dong, and C.~Long, ``Distributed mpc-based energy scheduling for islanded multi-microgrid considering battery degradation and cyclic life deterioration,'' \emph{Applied Energy}, vol. 329, p. 120168, 2023.

\bibitem{Ji_2023}
J.~Ji, M.~Zhou, R.~Guo, J.~Tang, J.~Su, H.~Huang, N.~Sun, M.~S. Nazir, and Y.~Wang, ``A electric power optimal scheduling study of hybrid energy storage system integrated load prediction technology considering ageing mechanism,'' \emph{Renewable Energy}, vol. 215, p. 118985, 2023.

\bibitem{Wang_2022}
T.~Wang, H.~Liang, B.~He, H.~Hua, Y.~Qin, and J.~Cao, ``A novel consensus-based optimal control strategy for multi-microgrid systems with battery degradation consideration,'' \emph{CSEE Journal of Power and Energy Systems}, pp. 1--13, 2022.

\bibitem{Satish_2}
S.~Vedula, M.~M. Bijaieh, E.~O. Boateng, and O.~M. Anubi, ``Degradation aware predictive energy management strategy for ship power systems,'' in \emph{2021 IEEE Electric Ship Technologies Symposium (ESTS)}, 2021, pp. 1--5.

\bibitem{Khalil_Book}
H.~Khalil, \emph{Nonlinear Control}, ser. Always Learning.\hskip 1em plus 0.5em minus 0.4em\relax Pearson, 2014.

\bibitem{SONG2018433}
Z.~Song, J.~Li, J.~Hou, H.~Hofmann, M.~Ouyang, and J.~Du, ``The battery-supercapacitor hybrid energy storage system in electric vehicle applications: A case study,'' \emph{Energy}, vol. 154, pp. 433--441, 2018.

\bibitem{boyd_vandenberghe_2004}
S.~Boyd and L.~Vandenberghe, \emph{Convex Optimization}.\hskip 1em plus 0.5em minus 0.4em\relax Cambridge University Press, 2004.

\bibitem{Boyd_ADMM}
S.~Boyd, N.~Parikh, E.~Chu, B.~Peleato, and J.~Eckstein, ``Distributed optimization and statistical learning via the alternating direction method of multipliers,'' \emph{Foundations and Trends® in Machine Learning}, vol.~3, no.~1, pp. 1--122, 2011.

\end{thebibliography}

\end{document}